%% file: main-arxiv-2.tex
\documentclass[11pt]{article}  

\usepackage{booktabs}
\usepackage{fancyhdr}
\usepackage{pifont}
\usepackage{url}
\usepackage{amsfonts}
\usepackage{latexsym}
\usepackage{amssymb}
\usepackage{epsf}
\usepackage{mdwlist}
\usepackage{lineno}
\usepackage{mathrsfs}
\usepackage{amsmath}
\usepackage{graphicx}
\usepackage{flafter}
\usepackage{epic,eepic}
\usepackage{cases}
\usepackage{color,soul}
\usepackage{xcolor}
\usepackage{colortbl}
\usepackage{float}
\usepackage{dsfont}
\usepackage{amssymb}
\usepackage{multirow}
\usepackage{titlesec}
\usepackage{titletoc}
\usepackage[top=1in, bottom=1.5in, left=1in, right=1in]{geometry}
\usepackage[algo2e,ruled,linesnumbered,vlined]{algorithm2e}
\usepackage[usenames,dvipsnames]{pstricks}
\usepackage{epsfig}
\usepackage{pst-plot}
\usepackage{pst-node}
\usepackage{pst-math}
\usepackage{hyperref}
\usepackage{subfigure}

\usepackage[style=authoryear,backend=bibtex,firstinits=true,maxnames=3,url=false,doi=false,isbn=false]{biblatex}
\bibliography{gmg}

\input{mydefinitions}

\begin{document}

\title{Numerical Study of Geometric Multigrid Methods on CPU--GPU Heterogeneous Computers}
\author{
Chunsheng Feng\thanks{Hunan Key Laboratory for Computation \& Simulation in Science \& Engineering, Xiangtan University, China. These authors were partially supported by  NSFC Project (Grant No. 91130002 and 11171281), Program for Changjiang Scholars and Innovative Research Team in University of China (No. IRT1179) and the Key Project of
Scientific Research Fund of Hunan Provincial Science and Technology Department (No. 2011FJ2011) in China.},
Shi Shu\footnotemark[1], 
Jinchao Xu\thanks{Department of Mathematics, Pennsylvania State University, PA, USA. This author was partially supported by NSFC-91130011 and NSF DMS-1217142.},
Chen-Song Zhang\thanks{NCMIS and LSEC, Academy of Mathematics and System Sciences, Chinese Academy of Sciences, Beijing, China. This author was partially supported by NSFC-91130011.}
}


\maketitle

\begin{abstract}
The geometric multigrid method (GMG) is one of the most efficient solving techniques for discrete algebraic systems arising from elliptic partial differential equations. GMG utilizes a hierarchy of grids or discretizations and reduces the error at a number of frequencies simultaneously. Graphics processing units (GPUs) have recently burst onto the scientific computing scene as a technology that has yielded substantial performance and energy-efficiency improvements. A central challenge in implementing GMG on GPUs, though, is that computational work on coarse levels cannot fully utilize the capacity of a GPU. In this work, we perform numerical studies of GMG on CPU--GPU heterogeneous computers. Furthermore, we compare our implementation with an efficient CPU implementation of GMG and with the most popular fast Poisson solver, Fast Fourier Transform, in the cuFFT library developed by NVIDIA.
\end{abstract}


%
%
\input{introduction}
%
\input{Preliminary}
%
\input{gmg_algorithm}
%
\input{complexity_analysis}
%
\input{numerical_experiment}

%
\input{conclusion}
%
\input{acknowledgements}

\printbibliography

\end{document}

%% file: mydefinitions.tex
\newtheorem{theorem}{Theorem}[section]

\newtheorem{example}{Example}[section]
 \newtheorem{remark}[theorem]{Remark}
 
\newenvironment{itemizex}%
  {\begin{itemize}%
    \setlength{\topsep}{0pt}%
    \setlength{\itemsep}{0pt}%
    \setlength{\parskip}{0pt}}%
  {\end{itemize}}

%% file: introduction.tex

\section{Introduction}\label{sec:intro}

Simulation-based scientific discovery and engineering design demand extreme computing power and high-efficiency algorithms~\parencite{Petter1992EMM,Petter2003Schwarz,Hey2009,Kaushik2011,Keyes2011}. This demand is one of the main forces driving the pursuit of extreme-scale computer hardware and software during the last few decades. It has become increasingly important for algorithms to be well-suited to the emerging hardware architecture. In fact, the co-design of \emph{architectures}, \emph{algorithms}, and \emph{applications} is particulary important given that researchers are trying to achieve exascale ($10^{18}$ floating-point operations per second) computing. Although the question of what is the best computer architecture to achieve exascale or higher remains highly debatable, many researchers agree that hybrid architectures make sense due to energy-consumption constraints. There are already quite a few heterogeneous computing architectures available, including the Cell Broadband Engine Architecture (CBEA), Graphics Processing Units (GPUs), and Field Programmable Gate Arrays (FPGAs)~\parencite{Carpenter2009,Brodtkorb2010,Wolfe2012}.

A GPU is a symmetric multicore processor that can be accessed and controlled by CPUs. The Intel/AMD CPU accelerated by NVIDIA/AMD GPUs is probably the most commonly used heterogeneous high-performance computing (HPC) architecture at present. GPU-accelerated supercomputers feature in many of the top computing systems in the HPC Top500~\parencite{HPCTop5002012} and the Green500~\parencite{Green5002012}. Some ``old'' supercomputers, such as JAGAUR (now known as TITAN) of the Oak Ridge National Laboratory, are being redesigned in order to incorporate GPUs and thereby achieve better performance. GPUs have evolved from fixed-pipeline application-specific integrated circuits into highly programmable, versatile computing devices. Under conditions often met in scientific computing, modern GPUs surpass CPUs in computational power, data throughput, and computational efficiency per watt by almost one order of magnitude~\parencite{Buck2007,Chen2009a,Nickolls2010a}.

Not only are GPUs the key ingredient in many current and forthcoming  petaflop supercomputers, they also provide an affordable desktop supercomputing environment for everyday usage, with peak computational performance matching that of the most powerful supercomputers of only a decade ago. General-purpose graphics processing units (GPGPU), as a high-performance computational device are becoming increasingly popular. Today's NVIDIA Fermi GPU and the upcoming (expected in December 2012) Intel Many Integrated Core (MIC) architecture are the most promising co-processors with high energy-efficiency and computation power. The Intel Knight Corner MIC (50 cores) is capable of delivering 1 Teraflop operation in double precision per second, whereas the peak performance of Tesla 2090 is 665 Gigaflop operations in double precision per second. On the other hand, as GPUs have a high-volume graphics market, it is expected by many experts to have a price advantage over MIC, at least immediately after its launch.

Probably one of the most discussed features of the MIC architecture is  that it shares the  x86 instruction set such that users often assume that they do not need to change their existing codebase in order to migrate to MIC. However, this assumption is subject to argument as even if legacy code can easily be migrated, whether the application it is then used for is able to  achieve the desired  performance is  questionable. Achieving scalable scientific applications in the exascale era is our ultimate goal. Hence, software, more importantly algorithms, must adapt to unleash the power of the hardware. Unfortunately, none of the processors envisioned at present will relieve today's programmers from the hard work of preparing their applications. In fact, power constraints will actually cause us to use simpler processors at lower clock rates for the majority of our work. As an inevitable consequence, improvements in terms of  performance  will largely arise from more parallel algorithms and implementations.

In all likelihood, different applications will benefit from any given architecture in specific ways. Thus a one-size-fits-all solution will almost certainly not arise. In many numerical simulation applications, the most time-consuming aspect is usually the solution of large linear systems of equations. Often, as they are generated by  discretized  partial differential equations (PDEs), the corresponding coefficient matrices are very sparse. The Laplace operator (or Laplacian) occurs in many PDEs that describe physical phenomena such as heat diffusion, wave propagation and electrostatics, and gravitational potential. For many efficient methods for solving discrete problems arising from PDEs, a fast Poisson solver is  a key ingredient in achieving a high level of efficiency~\parencite{Xu2010}. Numerical schemes based on fast Poisson solvers have been successfully applied to many practical problems among which are computed tomography, power grid analysis, and quantum-chemical simulation~\parencite{Kostler2007,Sturmer2008,Shi2009,Yang2011}.

Because of its plausible linear complexity---i.e., the low computational cost of solving a linear system with $N$ unknowns is $O(N)$---the GMG method is one of the most popular Poisson solvers~\parencite{Hackbusch.Hackbusch.1985ys,Bramble.J1993,Briggs.BriggsHenson.2000kx,Trottenberg.U;Oosterlee.C;Schuller.A2001,Brandt2011}. Although the GMG's applicability is limited as it requires explicit information on the hierarchy of the discrete system, when it can be applied,  GMG is far more efficient than  its algebraic version, the algebraic multigrid (AMG) method~\parencite{Brandt.BrandtMcCormick.1982fk,Brandt.Brandt.1986vn,Ruge.RugeStuben.1987uq,Trottenberg.U;Oosterlee.C;Schuller.A2001}. Another popular choice is the direct solver based on the fast Fourier transform or the FFT~\parencite{Cooley1965} on tensor product grids. The computational cost of the FFT-based fast Poisson solver is $O(N \log N)$, and FFT can easily be called from highly optimized software libraries, such as FFTW~\parencite{FFTW05} and the Intel Math Kernel Library (MKL). These advantages make FFT an extremely appealing method~\parencite{Sturmer2008,Lord2008} when it is applicable.

It is well-known that heterogeneous architectures pose new programming difficulties compared to existing serial and parallel platforms~\parencite{Chamberlain.ChamberlainFranklin.2007fk,Brodtkorb2010}. In this paper, we investigate the performance of fast Poisson solution algorithms, more specifically, GMG and FFT, on modern hybrid computer environment accelerated with GPUs. Considerable effort has been devoted to developing efficient solvers on GPUs for linear systems arising from PDEs and other applications; see, for example, \cite{Bolz2003,Bell2008,Bell2009,Barrachina2009,Jeschke2009,Cao2010,Elble2010,Georgescu2010,Bell.BellDalton.2011kx,Heuveline.HeuvelineLukarski.2011ys,Heuveline.HeuvelineLukarski.2011vn,Knibbe2011}.

The main purpose of this paper is to consider the following important questions, all of which are central to understanding geometric multigrid methods on GPU architecture:
\begin{itemizex}
\item What challenges multilevel iterative methods pose in addition to standard sparse iterative solvers do? Is it possible to achieve a satisfactory speedup on GPUs for multigrid algorithms? 
\item How does the performance of multigrid algorithms on GPUs compare with their performance on state-of-the-art CPUs?
\item How much of the computational power of GPUs can be utilized for multigrid algorithms? How cost-effective are CPU--GPU systems?
\end{itemizex}
We will consider answers to these questions based on carefully designed numerical experiments described herein. We will also compare a GPU-implementation of GMG with the optimized implementation of direct solvers based on FFT in double precision as a numerical experiment.

The rest of the paper is organized as follows: In Section~\ref{sec:pre}, we introduce the preliminary features of the hardware and algorithms under investigation. In Section~\ref{sec:gmg}, we give details about our implementation of GMG in a CPU--GPU heterogeneous computing environment. In Section~\ref{sec:complexity}, we analyze the complexity of the GMG algorithm. We report our numerical tests and analysis in Section~\ref{sec:numer}. We then summarize the paper with some concluding remarks in Section~\ref{sec:conclusion}. 

%% file: Preliminary.tex

\section{Preliminaries}\label{sec:pre}

Graphics processing units (GPUs) recently burst onto the scientific computing scene as an innovative technology that has demonstrated substantial performance and energy-efficiency improvements for many scientific applications. 

\subsection{A Brief Glance at GPU and CUDA}

\begin{figure}[htbp] 
   \centering
   \includegraphics[width=3in]{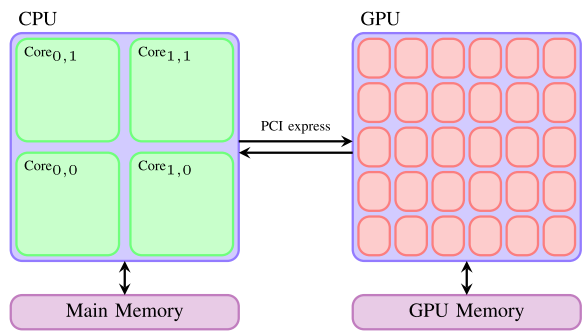}
   \caption{Schematic of heterogenous architecture: a quad core CPU in combination with a GPU.}
   \label{fig:memory}
\end{figure}

A typical CPU--GPU heterogenous architecture contains one or more CPUs (host) and a GPU (device). GPU has its own device memory, which is connected to the host via a PCI express bus. One of the main drawbacks of using such an architecture for PDE applications is that it is necessary to exchange data between the host and the device frequently; see Figure~\ref{fig:memory} \parencite{Brodtkorb2010}. Data must be moved to the GPU memory, and parallel kernels are launched asynchronously on the GPU by the host.

In regard to sparse matrix operations, the memory bandwidth is usually where the bottleneck occurs. The gap between the speed of floating-point operations and the speed for accessing memory grows every year~\parencite{Asanovic.AsanovicBodik.2006zr}. In this sense, we do not expect that the iterative linear solvers, which are usually memory-bounded, to readily derive benefits easily from increasing the number of cores. One way to address this problem is to use a high-bandwidth memory, such as Convey's Scatter-Gather memory~\parencite{Bakos2010}. Another is to add multithreading, where the execution unit saves the state of two or more threads, and can swap execution between threads in a single cycle---either by swapping between threads at a cache miss or by alternating between threads on every cycle. While one thread is waiting for memory, the execution unit keeps busy by switching to a different thread. 

CUDA (Compute Unified Device Architecture)~\parencite{CUDA4.1} is a parallel computing platform and programming model invented by NVIDIA. It delivers dramatic increases in computing performance by harnessing the power of the graphics processing unit (GPU). NVIDIA provides a complete toolkit for programming on the CUDA architecture, supporting standard computing languages such as C/C++ and Fortran. CUDA C and Fortran are the most widely used programming languages for GPU programming today~\parencite{Wolfe2012}. CUDA was developed simultaneously with the GeForce 8 architecture (NVIDIA's internal code name for the latter is Tesla), and publicly announced in  2006. In addition to CUDA, other options are OpenCL, AMD Stream SDK, and OpenACC supported by CAPS, CRAY, NVIDIA, and PGI.

\subsection{The Poisson Equation and Its Finite Difference Discretizations}
   Consider the Poisson equation
\begin{equation} \label{equ1-1-1}
\left\{
\begin{array}{rcl}
    - \Delta u &=&  f  \qquad \mbox{in }~\Omega \\
      u &=& 0  \qquad \mbox{on }\partial\Omega, \\
\end{array}
\right.
\end{equation}
where $\Omega = (0,1)^d \subset \mathbb{R}^d$.
The main reason why we choose this simplest possible setting is to emphasize that, even for a simple problem, the new heterogeneous architectures present challenges for numerical implementation. Another reason is to allow us to use explicit stencils and to avoid the bottleneck of sparse matrix-vector production~\parencite{Guo2012}. Furthermore, in this simple setting, we can compare our program with an existing efficient FFT implement in CUDA.  

The standard central finite difference method is applied to discretize (\ref{equ1-1-1})~\parencite{MR2153063}. In other words, the Laplace operator is discretized by the classical second-order central difference scheme. After discretization, we end up with a system of linear equations:
     \begin{equation} \label{linear1}
        \mathbf{A}\vec{u} = \vec{f}.
    \end{equation}

We use the five-point central difference scheme and the seven-point central difference scheme in 2D and 3D, respectively. Consider a uniform square mesh of $\Omega = [0, 1]^2$
with size $h = \frac{1}{n}$ and in which $x_i = ih, \, y_j = jh \, (i, j = 0,
1,\ldots, n)$. Let $u_{i,j}$ be the numerical approximation of $u(x_i , y_j )$. The five-point central
difference scheme for solving (\ref{equ1-1-1}) in 2D can be written as follows:
$$
{-u_{i-1,j}-u_{i,j-1}+4u_{i,j}-u_{i+1,j}-u_{i,j+1}}={h^2}f(x_i,y_j)
\qquad i, j=1,2,\ldots,n-1.
$$
Similar to the 2D case, we consider a uniform cube mesh of $\Omega = [0,
1]^3$ with size $h =\frac{1}{n}$ and in which $x_i
= ih$, $y_j = jh$ and $z_k = kh, ~i, j, k = 0, 1,\ldots , n$. Let
$u_{i,j,k}\approx u(x_i , y_j , z_k )$ be the approximate solution. The seven-point central difference scheme for solving (\ref{equ1-1-1}) in 3D reads
$$
 {-u_{i-1,j,k}- u_{i,j-1,k}- u_{i,j,k-1} + 6u_{i,j,k}- u_{i+1,j,k}- u_{i,j+1,k}- u_{i,j,k+1}}={h^2} f (x_i , y_j , z_k ),
$$
for all $i,j,k =1,2,\ldots,n-1$.

\subsection{Fast Fourier Transform}

A fast Fourier transform (FFT) is an efficient algorithm for computing the discrete Fourier transform (DFT) and its inverse. DFT decomposes a sequence of values into components of different frequencies. Computing DFT directly from its definition is usually  too slow to be practical. The  FFT can be used to compute the same result, but much more quickly. In fact, computing a DFT of $N$ points directly, according to its definition, takes $O(N^2)$ arithmetical operations, whereas FFT can compute the same result in $O(N \log N)$ operations~\parencite{Walker.Walker.1996ly}.

On tensor product grids, FFT can be used to solve the Poisson equation efficiently. We now explain the key steps for using FFT to solve the 2D Poisson equation (the 3D case is similar):
\begin{enumerate}
  \item Apply 2D forward FFT to $f(x,y)$ to obtain $\hat{f}(k_x,k_y)$, where $k_x$ and $k_y$ are the wave numbers. The 2D Poisson equation in the Fourier space can then be written as
\begin{equation}\label{fft4poisson}
    -\Delta u(x,y) = f(x,y) ~  \underrightarrow{~~FFT~~}~ -(k^2_x+ k^2_y  )\hat{u}(k_x,k_y) = \hat{f}(k_x,k_y).
\end{equation}
  \item Apply the inverse of the Laplace operator to $\hat{f}(k_x,k_y)$ to obtain $\hat{u}(k_x,k_y)$, which is the element-wise division in the Fourier space
      $$
       \hat{u}(k_x,k_y) = -\frac{\hat{f}(k_x,k_y)}{k^2_x+ k^2_y }.
      $$
  \item Apply 2D inverse FFT to $\hat{u}(k_x,k_y)$ to obtain $u(x,y)$.
\end{enumerate}

The NVIDIA CUDA Fast Fourier Transform (cuFFT version 4.1) library provides a simple interface for computing FFTs up to 10 times faster than MKL 10.2.3 for single precision.\footnote{cuFFT 4.1 on Telsa M2090, ECC on, MKL 10.2.3, and TYAN FT72-B7015 Xeon x5680 Six-Core 3.33GHz.} By using hundreds of processor cores on NVIDIA GPUs, cuFFT is able to deliver the floating point performance of a GPU without necessitating the development of custom GPU FFT implementation~\parencite{CUFFT}.

%% file: gmg_algorithm.tex

\section{Geometric Multigrid Method for GPU}\label{sec:gmg}

Multigrid (MG) methods in numerical analysis comprise a group of algorithms for solving differential equations using a hierarchy of discretizations. The main idea driving multigrid methodology is that of accelerating the convergence of a simple (but usually slow) iterative method by global correction from time to time, accomplished by solving corresponding coarse-level problems. Multigrid methods are typically applied to numerically solving discretized partial differential equations~\parencite{Hackbusch.Hackbusch.1985ys,Trottenberg.U;Oosterlee.C;Schuller.A2001}. In this section, we briefly review standard multigrid and full multigrid (V-cycle) algorithms and their respective implementations in a CPU--GPU heterogenous computing environment.

\subsection{Geometric multigrid method}

\begin{figure}[h!!]
\begin{center}
  \includegraphics[width=0.6\textwidth]{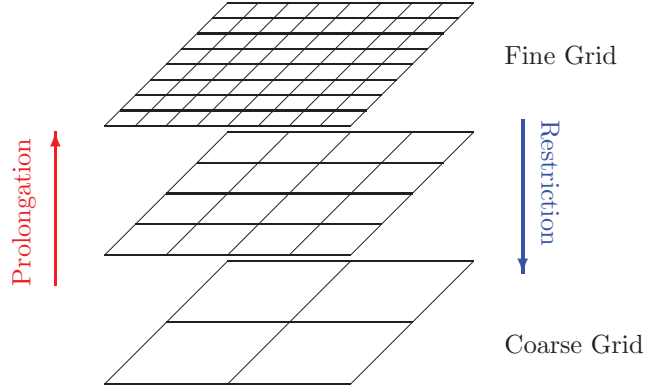}\\
%
%
%
\caption{Pictorial representation of a multigrid method with three grid levels.}\label{fig_multigrid}
\end{center}
\end{figure}

\begin{figure}[h!!]
\begin{center}
  \includegraphics[width=0.5\textwidth]{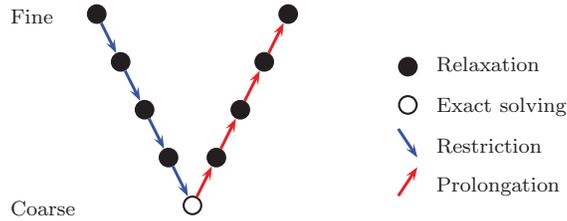}\\
  \caption{A schematic description of the V-cycle.}\label{fig_vcycle}
\end{center}
\end{figure}

The key steps in the multigrid method (see Figure~\ref{fig_multigrid}) are as follows:
\begin{itemizex}
  \item \textbf{Relaxation} or \textbf{Smoothing}: Reduce high-frequency errors using one or more smoothing steps based on a simple iterative method, like Jacobi or Gauss-Seidel.
  \item \textbf{Restriction}: Restrict the residual on a finer grid to a coarser grid.
  \item \textbf{Prolongation}: Represent the correction computed on a coarser grid to a finer grid.
\end{itemizex}
One of the simplest multilevel iterative methods is the multigrid V-cycle (see Figure~\ref{fig_vcycle}), in which the algorithm proceeds from top (finest grid) to bottom (coarsest grid) and back up again.

Suppose we have $L$ levels of nested spaces, numbered from $0$ to $L-1$. Let $\mathbf{R}^{l}_{l+1}$ and  $\mathbf{P}^{l+1}_{l}$ be the restriction (from level $l$ to $l+1$) and prolongation (from level $l+1$ to $l$) operators, respectively. The V-cycle algorithm (Figure 3) can be written as Algorithm~\ref{alg_V-cycle}. Note that, instead of using the standard recursive formulation, we use for loops, which coincides the actual implementation.

{\scriptsize
\begin{algorithm2e}[H]
\caption{Multigrid V-cycle $\vec{u} $ = MG-V($\mu_f,\mu_b, L, \vec{u}, \vec{f} $ ) }  \label{alg_V-cycle} 
\BlankLine
 \For{$l = 0 $ \KwTo $L -2$  }{

     Relax$_{forward}$( $\mu_f, \mathbf{A}_{l}, \vec{f}_{l}, \vec{u}_{l} $ )

     $\vec{r}_l = \vec{f}_l - \mathbf{A}_l \vec{u}_{l}$; $ \vec{f}_{l+1} = \mathbf{R}^{l}_{l+1}\vec{r}_l$
      }

     Relax$_{forward}$( $\mu_f, \mathbf{A}_{L-1}, \vec{f}_{L-1}, \vec{u}_{L-1} $ )

 \For{$l = L -1$ \KwTo $0$  }{

     $ \vec{u}_{l} = \vec{u}_{l} +\mathbf{P}^{l+1}_{l} \vec{u}_{l+1}$ 

     Relax$_{backward}$( $\mu_b, \mathbf{A}_{l}, \vec{f}_{l}, \vec{u}_{l} $ )
}
\end{algorithm2e}
}

\begin{remark}[Coarsest-level solver]\rm
Note that, for simplicity, we assume that the coarsest level $L-1$ contains one degree of freedom. Hence, Relax$_{forward}$( $\mu_f, \mathbf{A}_{L-1}, \vec{f}_{L-1}, \vec{u}_{L-1} $ ) in Algorithm~\ref{alg_V-cycle} solves the coarsest-level problem exactly. The same thing  happens in Algorithm~\ref{alg_fmg}.
\end{remark}

\begin{figure}[H]
\begin{center}
  \includegraphics[width=0.85\textwidth]{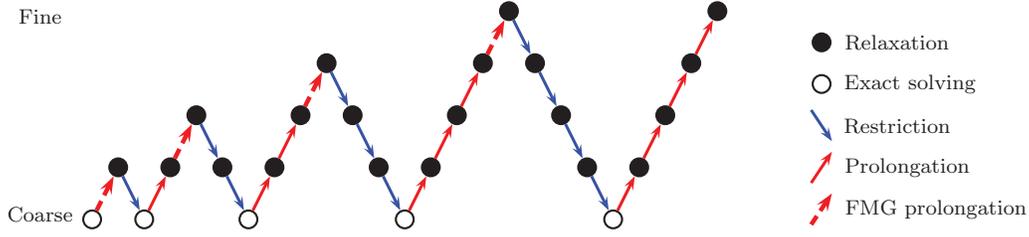}\\
  \caption{A schematic description of the full multigrid algorithm. The algorithm proceeds from left to right and from top (finest
grid) to bottom (coarsest grid).}\label{fig_fmg}
\end{center}
\end{figure}

The full multigrid (FMG) usually gives the best performance in terms of computational complexity. The idea of FMG is represented in Figure~\ref{fig_fmg}. Namely, we start from the coarsest grid and solve the discrete problem on the coarsest grid. Then, we interpolate this solution to the second-coarsest grid and perform one V-cycle. These two steps are repeated recursively on finer and finer grids, until the finest grid possible is achieved. The details are described in Algorithm~\ref{alg_fmg}.

{\scriptsize
\begin{algorithm2e}[H]
\caption{Full Multigrid V-cycle $\vec{u} $ = FMG-V($\mu_f,\mu_b, L,\vec{u}, \vec{f} $ ) }  \label{alg_fmg} 
\BlankLine

  init $\vec{f}_l,\vec{u}_l, \quad l = 0,\ldots, L-1$

  RelaxGS$_{forward}$( $\mu_f, \mathbf{A}_{L-1}, \vec{f}_{L-1}, \vec{u}_{L-1} $ )

 \For{$l = L-2 $ \KwTo $0$  }{

      $ \vec{u}_{l} =  \mathbf{P}^{l+1}_{l} \vec{u}_{l+1}$ 

      $u_l$ =V-cycle($\mu_f,\mu_b,l,\vec{u}_l, \vec{f}_l $ )
      }
\end{algorithm2e}
}

\subsection{Implementation of GMG on a CPU--GPU machine}

{\scriptsize
\begin{algorithm2e}[H]
\caption{ $\vec{u}, iter $ =  GMGSolve($d, L, L_{\theta},tol, maxit,\mu_f,\mu_b$) }  \label{GMGSolve_GPU}

\BlankLine

  \fboxrule=1.0pt \fcolorbox{gray}{gray!40}{\makebox[75mm] { for $l = 0$ \KwTo $L_{\theta}$ do: init $\vec{u}_l$ and $\vec{f}_l$; $resinit  = \|\vec{f}_{0} - \mathbf{A}\vec{u}_{0}\|  $ }}

   \fboxrule=1.0pt \fcolorbox{white}{white!40}{\makebox[78mm] { for $l = L_{\theta}$ \KwTo $L -1$ do: init  $\vec{u}_l$ and $\vec{f}_l$ ; res = resinit; iter = 0 }}


\While {$ (  res > tol \times resinit  )$  {\bf and} $( iter < maxit )$} {
  \For{$l = 0$ \KwTo $L_{\theta}-1$ on GPU}{
         \fboxrule=1.0pt \fcolorbox{gray}{gray!40}{\makebox[87mm] {   RelaxGS$_{forward}$( $\mu_f, \mathbf{A}_{l}, \vec{f}_{l}, \vec{u}_{l} $ ); $\vec{r}_l = \vec{f}_l - \mathbf{A}_l \vec{u}_{l}$; $\vec{f}_{l+1} = \mathbf{R}^{l}_{l+1}\vec{r}_l$ }}
        }

 copy $\vec{f}_{L_{\theta}}$ from DEVICE memory to HOST memory

 $\vec{u}_{L_{\theta}} $ = MG-V($\mu_f,\mu_b, L_{\theta}, \vec{u}_{L_{\theta}}, \vec{f}_{L_{\theta}} $ )  ~~~~~~~~~~~~~~~~~~~

 copy $\vec{u}_{L_{\theta}}$ from HOST memory to DEVICE memory

\For{$l = L_{\theta}-1$ \KwTo $ 0$ on GPU }{

          \fboxrule=1.0pt \fcolorbox{gray}{gray!40}{\makebox[87mm] {$ \vec{u}_{l} = \vec{u}_{l} + \mathbf{P}^{l+1}_{l} \vec{u}_{l+1}$; RelaxGS$_{backward}$( $\mu_b, \mathbf{A}_{l},  \vec{f}_{l},  \vec{u}_{l} $ )  }} 

        }
    \fboxrule=1.0pt \fcolorbox{gray}{gray!40}{\makebox[30mm] {   $res =\|\vec{f} - \mathbf{A}\vec{u}_{0}\| $ }}

        iter = iter + 1
}

 copy $\vec{u}_0$ from DEVICE memory to HOST memory $\vec{u}$

\end{algorithm2e}
}

\begin{remark}\rm We offer some remarks about our implementation:
\begin{enumerate}
\item There have been many discussions on how to implement geometric multigrid methods efficiently on modern computer architectures; see, for example, \cite{Weiss_2001}. We use a four-color and an eight-color Gauss-Seidel smoother for RelaxGS in 2D and 3D, respectively. As we are considering structured grids and the coloring is easy to obtain, we prefer the GS smoother over the weighted Jacobi smoother (see Table~\ref{table_smoother} for numerical comparison of several different smoothers). A weighted Jacobi smoother  is likely to achieve a higher peak performance and higher speedup over the corresponding CPU version; however, it usually requires more iterations and wall time compared with the colored GS smoother if both methods use same multilevel iteration, like V-cycle.
\item The \emph{gray} boxes represent the code segments running on GPU (kernel functions). When ${L_\theta = 0}$, Algorithm~\ref{GMGSolve_GPU} runs on CPU completely, and when ${L_\theta = L}$, Algorithm~\ref{GMGSolve_GPU} runs on GPU solely. However, when ${0 < L_\theta < L}$, these functions run in CPU--GPU hybrid mode. Note that we prefer this hybrid algorithm over using a U-cycle with a direct solver on GPU---This is because the matrix-free geometric multigrid is much faster than the direct solver; Furthermore, to use a direct solver on GPU, we need to form the coarse level matrices, which will cost extra time. 
\item Graphics processors provide texture memory to accelerate frequently performed operations. As optimized data access is crucial to GPU performance, the use of texture memory can sometimes provide a considerable performance boost; see Table~\ref{table_texture} for numerical results. We band the vectors $\vec{u}_l$ as one dimension texture memory in function $\vec{r}_l = \vec{f}_l - \mathbf{A}_l \vec{u}_{l} , ~ l= 0,1,\ldots,L-2$.
\end{enumerate}
\end{remark}



%% file: complexity_analysis.tex

\section{Complexity Analysis of the GMG Algorithm}\label{sec:complexity}


For the geometry multigrid of the finite difference method on structured meshes, it is not necessary to explicitly assemble the global stiffness matrix $A$ (i.e., matrix-free).
A subroutine for the matrix-vector multiplication of the corresponding finite difference operator is called whenever we need to compute $\vec{r} =\vec{f} - \mathbf{A}\vec{u}$. This subroutine can be implemented directly from the central difference scheme. One V-cycle of the GMG algorithm consists of computing the residual, forward relaxation, backward relaxation, restriction, prolongation, and the inner product. Take the 2D case as an example: we can analyze the time and space complexity of these operations. The time complexity of one V-cycle can be proven to be $O(N)$ for both the 2D and 3D cases. The optimal complexity of GMG has been analyzed by Griebel~\parencite{Griebel.Griebel.1989fk}. Next, we give the exact operation counts for each component in a V-cycle for the convenience of later discussions in \S\ref{sec:numer}. To the best of our knowledge, we cannot find an appropriate reference which can serve for this purpose. 

\subsection{2D case}

(1) Computing residual $\vec {r}_l =\vec{f}_l - \mathbf{A}_l\vec {u}_l$:
      \begin{eqnarray}
        r^l_{i,j}  = f^l_{i,j}-4u^l_{i,j} + u^l_{i \pm 1,j} + u^l_{i,j\pm1}, \quad i,j = 1,2,\ldots, n_l-1, \label{residual_2D}
      \end{eqnarray}
   where $u^l_{i \pm 1,j} =u^l_{i - 1,j} + u^l_{i+1,j}$ and in this case when $i - 1 = 0 $ or $i + 1 = n_l +1 $, then $u^l_{i \pm 1,j} = 0$. From now on, we will use the notation $\pm$  in this section. The equation (\ref{residual_2D}) requires 6 floating-point operations or work units (W) per unknown in the 2D case. Hence, the  total number of floating-point operations for the residual in one V-cycle is
$      
        W_{Residual}  = 6\sum_{l=0}^{L-2} (n_l)^{2} + 6n_0^2.     \label{total_residual_2D}
$

\smallskip\noindent
(2) Pre- and post-smoothing using 4-color Gauss-Seidel relaxations:
      \begin{eqnarray}
         u^l_{i,j}  = \frac{1}{4}(f^l_{i,j} + u^l_{i \pm 1,j} + u^l_{i,j\pm1}), \quad i,j = 1,2,\ldots, n_l-1 . \label{relax_2D}
      \end{eqnarray}
Equation (\ref{relax_2D}) shows  that there are 5 floating-point operations per unknown in the 2D case and the total number of floating-point operations in the forward and backward Gauss-Seidel relaxation is
$      
         W_{GSforward}  = 
         W_{GSbackward}  = 5\sum_{l=0}^{L-2}(n_l)^{2}.   \label{btrelax_2D}
$%

\smallskip\noindent
(3) Restriction operator $\vec {r}_{l+1} = \mathbf{R}^{l}_{l+1} \vec {r}_l$:
\begin{eqnarray}
      r^{l+1}_{i,j}  = \frac{1}{8}(2r^l_{2i,2j } + r^l_{2i \pm 1,2j} + r^l_{2i,2j\pm1}+ r^l_{2i - 1,2j- 1}+r^l_{2i+ 1,2j +1}), \quad i,j = 1,2,\ldots, n_l-1. \label{Restriction2d} 
      \end{eqnarray}
Equation (\ref{Restriction2d}) requires 8 floating-point operations per unknown in the 2D case. Furthermore, we obtain  the total floating-point operations of restriction for one V-cycle as
$      
        W_{Resitriction}  = 8\sum_{l=1}^{L-1}(n_l)^{2}.    \label{tRestriction2d} 
$

\smallskip\noindent
(4) Prolongation operator $\vec {e}_{l} =\vec {e}_{l}+ \mathbf{P}^{l+1}_{l} \vec {e}_{l+1}$:
  \begin{eqnarray}
  \label{Prolongation2d}    
   \begin{array}{ll}
      e^{l}_{2i,2j} = e^{l}_{2i,2j} + e^{l+1}_{i,j}, &  e^{l}_{2i+1,2j} = e^{l}_{2i+1,2j} + \frac{1}{2}(e^{l+1}_{i,j}+e^{l+1}_{i+1,j}),\\
       e^{l}_{2i,2j+1} = e^{l}_{2i,2j+1} + \frac{1}{2}(e^{l+1}_{i,j}+e^{l+1}_{i,j+1}), &
       e^{l}_{2i+1,2j+1} = e^{l}_{2i+1,2j+1} + \frac{1}{2}(e^{l+1}_{i,j}+e^{l+1}_{i+1,j+1}).\\
       ~~&  i,j =1,\ldots, n_{l+1}-1
    \end{array} 
  \end{eqnarray}
Equation (\ref{Prolongation2d}) shows that there are $\frac{3 \times 3 + 1}{4}$ floating-point operations per unknown for the 2D case. Furthermore, we can obtain that the total number of floating-point operations of prolongation for one V-cycle is
$
       W_{Prolongation}  = 2.5 \sum_{l=0}^{L-2}(n_l)^{2}.   \label{tProlongation2d}
$

\smallskip\noindent
(5) Computing the norm of the residual:
      \begin{eqnarray}
           \| \vec {r}_0 \|_{L^2} = \sum\limits_{j=1}^{(n_0)^2}r^0_{j}r^0_{j}. \label{Inner_product2D} 
      \end{eqnarray}
We can see that the total number of floating-point operations for computing $\ell^2$-norm is $2(n_0)^{2}$.

By combining the above five components, we can get the total number of floating-point operations per unknown for the 2D case in one V-cycle as
  \begin{eqnarray*}
    \frac{ 6\sum\limits_{l=0}^{L-2}(n_l)^{2}
    + 5\sum\limits_{l=0}^{L-1}(n_l)^{2}
    + 5\sum\limits_{l=0}^{L-2}(n_l)^{2}
    + 8\sum\limits_{l=1}^{L-1}(n_l)^{2}
    +  \frac{10}{4}\sum\limits_{l=0}^{L-2}(n_l)^{2}
    + (2+6)(n_0)^{2}  }    {(n_0)^2} \cong 36.
  \end{eqnarray*}
This means the total number of floating-point operations per unknown required by one V-cycle in  the 2D case is about $36$. 

\subsection{3D case}

Similarly, we can count the complexity of a V-cycle in 3D:
    \begin{itemize}
      \item[(1)] Computing residual: $W_{Residual}  = 8\sum_{l=0}^{L-2}(n_l)^{3} + 8 (n_0)^{3}$.

      \item[(2)] Smoothing:
$W_{GSforward} = W_{GSbackward}  = 7\sum_{l=0}^{L-2}(n_l)^{3}.$
     
      \item[(3)] Restriction operator: $W_{Resitriction}=16\sum_{l=1}^{L-1}(n_l)^{3}.$

      \item[(4)] Prolongation operator: $W_{Prolongation}  = \frac{23}{8} \sum_{l=0}^{L-2}(n_l)^{3}.$

       \item[(5)] Computing the norm of the residual: $W_{Norm} = (n_0)^3$.
       
    \end{itemize}
By combining the above estimates, we obtain the total number of floating-point operations per unknown for the 3D case in one V-cycle as
  \begin{eqnarray*}
    \frac{ 8\sum\limits_{l=0}^{L-2}(n_l)^{3}
    + 7\sum\limits_{l=0}^{L-1}(n_l)^{3}
    + 7\sum\limits_{l=0}^{L-2}(n_l)^{3}
    + 16\sum\limits_{l=1}^{L-1}(n_l)^{3}
    +  \frac{22}{8}\sum\limits_{l=0}^{L-2}(n_l)^{3}
    + (2+8)(n_0)^3  }    {(n_0)^3} \cong 41.
  \end{eqnarray*}
Hence, the total number of floating-point operations per unknown for one V-cycle is $36$ and $41$ for the 2D and 3D cases, respectively.

\begin{remark}[Space complexity of V-cycle in GMG]\rm
In Algorithm \ref{GMGSolve_GPU}, we need only keep $\vec{u}_l$, $\vec{f}_l \, (l=0,1, \ldots, L-1)$  and $\vec{r}_l \, (l=0,1,\ldots,L-2)$
  in the host or device memory. Therefore, we obtain the memory space complexity of GMG (Algorithm~\ref{GMGSolve_GPU}) as follows:
  \begin{equation}
   \text{Memory}/N
    = \frac{1}{(n_0)^d} \left\{2 \sum\limits_{l=0}^{L-1}(n_l)^{d}
       + \sum\limits_{l=0}^{L-2}(n_l)^{d} \right\}   \cong  4. \label{memfor2D}
  \end{equation}
Equation (\ref{memfor2D}) shows that the memory space complexity of GMG Algorithm  \ref{GMGSolve_GPU}  has about  4 times as many unknowns in both 2D and 3D. In fact, in our numerical experiments, we find the memory usage is about $3.7N$ in 2D and $3.3N$ in 3D. 
\end{remark}

%% file: numerical_experiment.tex

\section{Numerical Experiment}\label{sec:numer}

 In this section, we perform several numerical experiments and analyze the performance of GMG as proposed in Algorithm~\ref{GMGSolve_GPU}. In order to obtain relatively accurate wall times, all reported computing times are the averages of 100 runs. In order to provide a fair comparison, we perform experiments to compare our implementation of GMG in CUDA with an efficient OpenMP version of GMG and a direct solver based on cuFFT. In order to eliminate effects of implementation as much as possible, we consider the simplest test problem, namely the Poisson equation on uniform grids:
\begin{example}[Poisson Equation] \label{example3D}
In the model problem \eqref{equ1-1-1} for $d=2,3$, we take the right-hand side
$$ f(x)= \Pi_{i=1}^d\sin(\pi x_i), \qquad x \in \Omega=(0,1)^d \subset \mathbb{R}^d.$$
The tolerance for the convergence of Algorithm~\ref{GMGSolve_GPU} is $\mathrm{tol} = 10^{-6}$.
\end{example}
%

\subsection{Environment for Comparisons}

Our focal computing environment is a HP workstation with a low-cost commodity-level NVIDIA GPU. Details in regard to the machine are set out in Table~\ref{CentOS6.2}. 
For numerical experiments, we use an AMD FX(tm)-8150 Eight-Core 3.6GHz CPU (its peak performance in double precision is 1.78GFLOPs\footnote{Obtained experimentally using LINPACK (\url{http://www.netlib.org/benchmark/linpackc}). }) and an NVIDIA GeForce GTX~480 GPU. GTX~480 supports CUDA and it is composed of 15 multiprocessors, each of which has 32 cores (480 cores in total). Each multiprocessor is equipped with 48KB of very fast shared memory, which stores both data and instructions. All the multiprocessors are connected to the global memory, which is understood as an SMP architecture. The global memory is limited to a maximum size of 1.5GB. 
However,
there is also a read-only cache memory called a texture cache, which is bound to a part of the global memory when a code is initiated by the multiprocessors. 
\begin{table}[H] 
  \centering
  \small
  \caption{Experiment Environment}\label{CentOS6.2}
  \begin{tabular}{ll} \hline
          CPU Type      & AMD FX-8150 8-core\\
          CPU Clock  & 3.6 GHz $\times$ 8 cores\\
          CPU Energy Consumption & 85 Watts (idle) $\sim$ 262 Watts (peak) \\
          CPU Price & 300 US Dollars\\
          Host Memory Size & 16GB \\
          \hline
          GPU Type     & NVIDIA GeForce GTX~480sp\\
          GPU Clock & 1.4 GHz $\times 15 \times 32$ cores \\
          GPU Energy Consumption & 141 Watts (idle) $\sim$ 440 Watts (peak) \\
          GPU Price & 485 US Dollars \\
          Device Memory Size & 1.5GB \\
          \hline
          Operating System  & CentOS 6.2 \\
          CUDA Driver & CUDA 4.1 \\
      Host Compiler & gcc 4.4.6 \\
    Device Compiler & nvcc 4.1 \\ \hline
  \end{tabular}
\end{table}

\begin{remark}[Cost Effectiveness]\label{rem:cost}\rm
From Table~\ref{CentOS6.2}, we notice that the initial cost and the peak energy-consumption cost of this particular system with GPU is roughly 2 to 2.5 times of the system without GPU. These extra initial and energy costs must be considered in our comparisons. 
\end{remark}

\subsection{CPU v.s. GPU}

Before we embark on comparing the GMG algorithm on CPU and GPU, we collect a few benchmarks. The main performance parameters of GeForce GTX~480 is described  in Table \ref{Theoretical_peak_gtx480}\footnote{Numbers in the last four rows are obtained experimentally using the \emph{bandwidthtest} of CUDA 4.1 SDK.}. On this system, we test the performance of double precision matrix-vector products\footnote{Using explicit stencil to avoid storing sparse coefficient matrices.} and the results are reported in Table~\ref{table_residual_GFLOPs_2D} and \ref{table_residual_GFLOPs_3D} for 2D and 3D, respectively\footnote{The CPU version is implemented using OpenMP to take advantage of the multicore platform.}. Here, since we do not store the coefficient matrix, the performance of matrix-vector production is much better than sparse matrix-vector products. Since (matrix-free) matrix-vector products are responsible to most of the computation work in the geometric multigrid, we wish that our implementation of GMG on GPUs can achieve this performance (roughly 30 GFLOPs in double precision).

\begin{table}[H]
\centering
   \caption{Theoretical peak performance of NVIDIA GeForce GTX~480}
  \label{Theoretical_peak_gtx480}
   \vskip 0.2cm
    \begin{tabular}{lr} \hline
   Double precision performance [GFLOPs]  &177.00  \\ 
   Theoretical memory bandwidth [GB/s]      &  177.00 \\ \hline
          Device to device  memory bandwidth [GB/s] & \qquad 148.39 \\ 
          Device to host  memory bandwidth [GB/s] & 4.46 \\ 
          Host to host  memory bandwidth [GB/s]     & 9.44   \\ 
          Host to device  memory bandwidth [GB/s] & 3.92 \\ \hline
    \end{tabular}
\end{table}

\begin{table}[H]
\centering\caption{Wall times, GFLOPs (double precision), and speedups for computing the residual  in 2D}\label{table_residual_GFLOPs_2D}
\begin{tabular}{c|cc|cc|c}\hline
\multirow{2}{*}{L}  &\multicolumn{2}{c|}{CPU}&\multicolumn{2}{c|}{GPU} & \multirow{2}{*}{Speedup}
 \\ \cline{2-5}
 & Wall Time (sec) & GFLOPs & Wall Time (sec) &GFLOPs &   \\\hline
    8      &1.607e$-$4   &2.47         &1.677e$-$5        &23.63         &9.57 \\ 
    9      &6.296e$-$4   &2.51         &5.027e$-$5        &31.41         &12.51 \\ 
    10     &3.079e$-$3   &2.05         &1.819e$-$4        &34.66         &16.91 \\ 
    11     &1.269e$-$2   &1.99         &7.061e$-$4        &35.68         &17.93  \\ 
    12     &5.107e$-$2   &1.97         &2.806e$-$3        &35.89         &18.22 \\ \hline
\end{tabular}
\end{table}
\begin{table}[H]
\centering\caption{Wall times, GFLOPs (double precision), and speedups for computing the residual  in 3D}\label{table_residual_GFLOPs_3D}
\begin{tabular}{c|cc|cc|c}\hline
\multirow{2}{*}{L}  &\multicolumn{2}{c|}{CPU}&\multicolumn{2}{c|}{GPU} & \multirow{2}{*}{Speedup}
 \\ \cline{2-5}
 & Wall Time (sec) & GFLOPs & Wall Time (sec) &GFLOPs &   \\\hline
    5     &7.111e$-$5   &4.04         &2.067e$-$5        &13.91         &3.44   \\ 
    6     &5.613e$-$4   &3.91         &7.293e$-$5        &30.13         &7.70  \\ 
    7     &5.565e$-$3   &3.09         &4.547e$-$4        &37.77         &12.22   \\ 
    8     &5.941e$-$2   &2.29         &4.536e$-$3        &29.94         &13.07  \\ \hline
\end{tabular}
\end{table}

First, we note that, with the given tolerance, both the GPU version and the CPU version of GMG achieve the optimal discretization error, $O(h^2)$, in $\mathbb{R}^2$ and $\mathbb{R}^3$. Furthermore,  the GPU and the CPU versions take the same number of iterations as each other to reach the given convergence tolerance, i.e., the GPU version is equivalent to the corresponding serial version. We also notice that, in all our numerical comparisons of GMG on CPU and on GPU, we use only one sweep of four-color Gauss--Seidel (GS) relaxation as smoother (i.e. $\mu_f = \mu_b = 1$) because it yields both fast convergence rate and good parallelism scalability. For example, Table~\ref{table_smoother} shows a comparison of GMG with a few other widely-used smoothers on GPU, like the Weighted Jacobi (WJ) method. From this table, we can also see that the gap between the WJ method and the 4-color GS decreases as the problem size increases.
\begin{table}[H]
        \centering\caption{Comparison of four different smoothers in 2D on GPU ($L_{\theta}=L$). Here \#It is the number of iterations, and wall times are in seconds. }\label{table_smoother}
        \smallskip
        \begin{tabular}{c|cc|cc|cc|cc}\hline
     \multirow{2}{*} {$L$ } & \multicolumn{2}{c|}{WJ (weight=$0.667$)} & \multicolumn{2}{c|}{WJ (weight=$0.8$)} & \multicolumn{2}{c|} {2-Color GS}  & \multicolumn{2}{c} {4-Color GS}
        \\ \cline{2-9}
       &\#It & Wall Time  &\#It & Wall Time  &\#It & Wall Time &\#It & Wall Time   \\\hline

        8  &22 &4.84e$-$2 &18  &3.97e$-$2   &16  &7.02e$-$3  &11  &6.68e$-$3           \\ 

        9  &22 &7.10e$-$2 &18  &5.89e$-$2   &16  &1.85e$-$2  &11   &1.46e$-$2    \\ 

        10 &22 &1.20e$-$1 &19  &1.01e$-$1   &16  &5.15e$-$2  &11   &3.74e$-$2   \\ 

        11 &22 &3.04e$-$1 &19  &2.66e$-$1   &16  &1.79e$-$1  &11   &1.28e$-$1    \\ 

        12 &22 &1.02e$+$0 &19  &8.88e$-$1   &16  &6.82e$-$1  &11   &4.85e$-$1    \\ \hline

       \end{tabular}
\end{table}

\begin{figure}[H]
  \centering
  \includegraphics[width=0.99\textwidth]{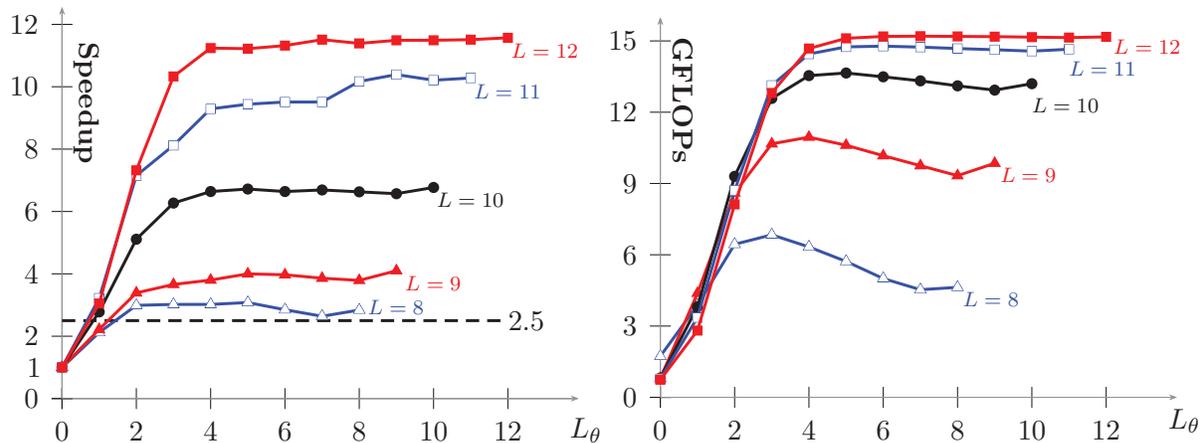}\\
\caption{Speedup (left) of GMG on GPU compared with its CPU OpenMP version and performance (right) of GMG on GPU in 2D.}\label{figure_2D_GMG}
\end{figure}

\begin{figure}[H]
  \centering
  \includegraphics[width=0.99\textwidth]{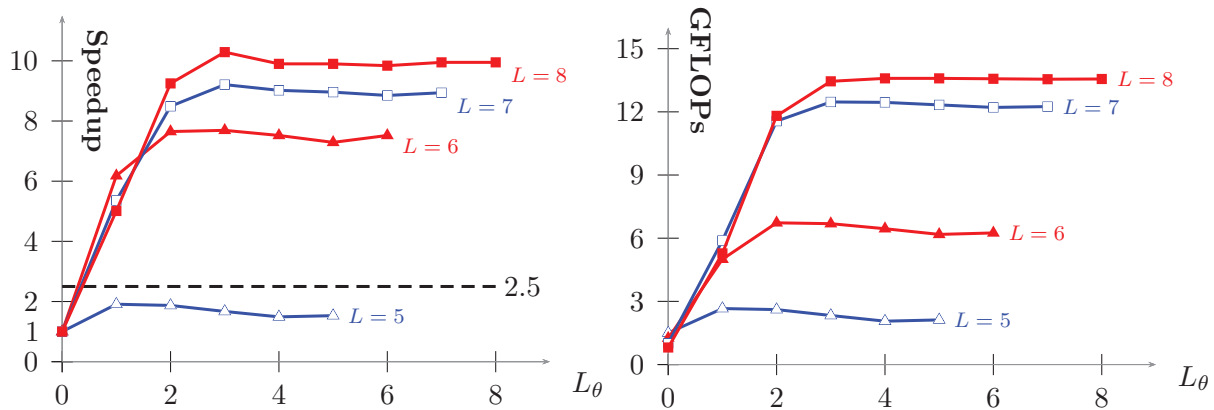}\\
\caption{Speedup (left) of GMG on GPU compared with its CPU OpenMP version and performance (right) of GMG on GPU in 3D.}\label{figure_3D_GMG}
\end{figure}

Second, we compare the total wall times for the GPU version and the CPU version (with OpenMP) of GMG. For the 2D test problem, in the best-case scenario, the GPU version of GMG can achieve about $11.5$ times speedup (see Figure~\ref{figure_2D_GMG}, Left). Moreover, the speedup increases as $L$ increases, but with a plateau zone when $L_{\theta} \ge 4$. This indicates that, when problem size is large, it is more efficient to carry out all computational work on GPU instead of sending data back and forth between GPU and CPU. For $L=12$, the implementation achieves $15.2$ GFLOPs (Figure~\ref{figure_2D_GMG}, Right) in double precision, which is $8.6\%$ of the theoretical peak performance of GTX~480 or $42\%$ of the performance of (finest grid) matrix-vector product on GPU (see Table~\ref{table_residual_GFLOPs_2D}). Similarly, for the 3D test problem, in the best-case scenario, we can achieve $10.3$ times speedup and 15 GFLOPs in double precision, which is $51\%$ of the matrix-vector product operation on the finest grid
(see Figure~\ref{figure_3D_GMG} and Table~\ref{table_residual_GFLOPs_3D}). Similar to the 2D case above, if the size of the problem is large enough, then we should run it completely on GPU. 


%
%

In Figure~\ref{figure_2D_GMG} (Left) and \ref{figure_3D_GMG} (Left), we draw two dashed lines to highlight $2.5\times$ speedup. As we have mentioned earlier in Remark~\ref{rem:cost}, the energy cost of our system with the GTX~480 GPU running at full power is roughly 2.5 times the system without using it. Hence, the speedup of the algorithm on GPU has to be at least above this dashed line to make it economically sound. For example, for $L=5$ in 3D, using GMG on GPU is not efficient enough to justify the energy consumption. On the other hand, one should also be aware of that, since the problem size is small for $L=5$ in 3D, many of the stream processors are idle during computation and the energy consumption of the GPU card might not be as high as the peak consumption indicated in Table~\ref{CentOS6.2}.

\begin{remark}[Multicore effect]\rm
We note that, in all our numerical tests on CPUs, we use all available CPU core. It is fair to say that on AMD FX8150 (8-core) the speedup of eight-thread GMG (over one thread version when compiler optimization O2 is enabled) is, in general, 2.0 to 3.0, depending on the algorithm and implementation.
\end{remark}

\subsection{Performance of GMG on GPUs}\label{ssc:gpu}

In this subsection, we would like to understand more about which part(s) of the GPU implementation are the bottleneck(s).
Tables \ref{table_ker_trans_2D} and \ref{table_ker_trans_3D} show the computing time and communication time  in 2D and 3D, respectively. In these two tables, total time = computing time + communication time, where the communication time refers to the time for transferring data between the host and the device and the computing time means the time for everything else. Here, since we take $L_\theta=L$, the computing time is approximately equal to the kernel time on GPU. The numerical results show that the CUDA kernel computation takes $75\%$ to $90\%$ of the total wall time. Furthermore, as one can expect, the portion of communication time decreases as problem size increases.

\begin{table}[H]
\centering \caption{Computing time (seconds) and communication time (seconds)  in 2D }\label{table_ker_trans_2D}
    \begin{tabular}{l|ccccc} \hline
       $L\; (L_{\theta} = L)$ &8             &9             &10               &11                 & 12 \\ \hline
    Computing time            & 5.934e$-$3     & 1.064e$-$2     & 3.134e$-$2        &1.124e$-$1           &4.339e$-$1       \\
    Communication time     & 9.070e$-$4     & 3.491e$-$3     & 5.884e$-$3        &1.522e$-$2           &5.169e$-$2        \\
    Communication/Total   & 13.26\%         & 24.70\%         & 15.81\%            &11.93\%               &10.64\%       \\ \hline
    \end{tabular}\\
\end{table}
\begin{table}[H]
\centering\caption{Computing time (seconds) and communication time (seconds)  in 3D  }\label{table_ker_trans_3D}
    \begin{tabular}{l|cccc} \hline
     $L\; (L_{\theta} = L)$      &5             &6             &7              &8   \\ \hline
    Computing time               &8.544e$-$3      &2.404e$-$2      &9.997e$-$2       &7.287e$-$1        \\
    Communication time        &4.938e$-$4      &3.716e$-$3      &9.574e$-$3       &5.141e$-$2        \\
    Communication/Total      &5.46\%          &13.39\%          &8.74\%          &6.59\%            \\ \hline
    \end{tabular}
\end{table}

Table~\ref{timerate_each_func} shows the wall time percentage (ratio to the total kernel time) for each function in one V-cycle. Table~\ref{gflops_each_func} shows the number of GFLOPs for each function of one V-cycle. By comparing the numerical results on CPU  ($L_\theta=0$) and GPU ($L_\theta=L$), we notice that the multicolored Gauss-Seidel smoother counts for more than $50\%$ of the total kernel time on GPU and it yields less speedups. Furthermore, because we are using the multicolored GS smoother, we need to launch the GS kernel several times (equal to the number of colors)---this introduces some overhead. On the other hand, as we pointed out earlier, using the weighted Jacobi method does not help because of the deteriorated convergence rate, although it results in better parallelism. Another observation is that computing the Euclidean norm of the residual gets very low efficiency on GPU due to it requires a summation of large amount of floating-point numbers.

\begin{table}[H]
\centering\caption{Floating-point operation complexity in terms of total degrees of freedom ($N$) and wall time percentage of each subroutine in one V-cycle. All computations are done on CPU if $L_\theta=0$ and done on GPU if $L_\theta=L$.}\label{timerate_each_func}
    \begin{tabular}{lll|rrrrr} \hline
              &   \multicolumn{2}{c|}{Setting}                &Residual             &Smoother                          &Restriction&Prolongation&Norm   \\ \hline
\multirow{3}{*}{2D} & \multicolumn{2}{c|}{Complexity}  & $14N$ & $13.4N$ & $2.67N$ & $3.34N$ & $2N$ \\
& $L=12$ & $L_{\theta} = 0$    &30.70\%     &45.02\%              &6.54\%  &15.84\% &1.91\%        \\ 
     & $L=12$ & $L_{\theta} = L$    &22.11\%     &50.40\%              &5.99\%  & 9.92\%&11.59\%         \\ \hline
\multirow{3}{*}{3D} & \multicolumn{2}{c|}{Complexity}  & $17.2N$ & $16.1N$ & $2.30N$ & $3.31N$ & $2N$ \\
& $L=8$ & $L_{\theta} = 0$     &28.60\%     &54.08\%              &3.90\%  &11.60\% &1.83\%        \\ 
     & $L=8$ & $L_{\theta} = L$     &25.34\%     &53.99\%              &4.38\%  &15.69\% &2.79\%        \\ \hline
    \end{tabular}
\end{table}

\begin{table}[H]
\centering\caption{Double precision performance (GFLOPs) of each subroutine in one V-cycle. All computations are done on CPU if $L_\theta=0$ and done on GPU if $L_\theta=L$.}\label{gflops_each_func}
    \begin{tabular}{lll|rrrrr} \hline
                 & Setting  &               &Residual             &Smoother                          &Restriction&Prolongation&Norm   \\ \hline
    \multirow{2}{*}{2D}  & $L=12$ & $L_{\theta} = 0$    &2.04    &1.30             &1.78  &0.91   &4.99        \\ 
     & $L=12$ & $L_{\theta} = L$    &27.12     &11.33          &19.09  &15.87 &7.39         \\ \hline
    \multirow{2}{*}{3D}  & $L=8$ & $L_{\theta} = 0$     &2.31     &1.10              &2.25  &1.06 &4.08        \\ 
     & $L=8$ & $L_{\theta} = L$     &24.81     &10.94              &25.78  &19.43 &7.66        \\ \hline
    \end{tabular}
\end{table}

In Section~\ref{sec:complexity}, we analyzed the operation counts of each component in a V-cycle of GMG. Number of memory accesses is roughly equal to the number of floating-point operations in these functions. However, in Table~\ref{timerate_each_func} and \ref{gflops_each_func}, we observe that, for the modern desktop environments (multicore CPUs or GPUs), operation and memory access counts alone can no longer provide good measurements of practical efficiency of the algorithm. For example, in 2D, computing residual and the GS relaxation require similar amount of floating-point operations and have almost identical memory access pattern; however, due to the multicolor GS smoother has less build-in parallelism and introduces bigger overhead, it costs much more time than computing residuals, especially on GPU. This also suggests that one would prefer a more scalable smoother (like weight Jacobi) than a more effective smoother (like multicolored GS) when problem size becomes very large\footnote{This trend has already been observed in Table~\ref{table_smoother}.}.

Moreover, visiting coarser levels in the hierarchy also poses an additional difficult for efficient implementation of multilevel iterative methods on GPUs. As an example, we show a typical profile of GMG (3D Poisson with 2M unknowns) in Table~\ref{tab:profile}. On the finest levels (0-th and 1st), computing residual and smoothing can achieve relatively high performance. But, when reaches coarse levels, the performance of GMG drops dramatically due to less active threads are used there. As we mentioned before, we note that replacing the coarse level V-cycle with a direct solver on GPU (like~\cite{magma}) is not an option as this will slow down the method even more. 
\begin{table}[H]
\small
\centering\caption{Performance (GFLOPs) of major components of V-cycle on each level in 3D ($L=7$).}\label{tab:profile}
\begin{tabular}{|c|c|c|c|c|c|} \hline
  Level &   Wall Time (s)  & Percentage (\%) & Overall & Residual & Smoother \\ \hline
  0 &4.88e-03 &74.60 &16.20 &36.94 &11.50\\
  1 &4.95e-04 & 7.57 &14.91 &26.26 &11.52\\
  2 &2.02e-04 & 3.09 &4.780 & 9.66 & 3.91\\
  3 &1.46e-04 & 2.24 &0.906 & 2.31 & 0.79\\
  4 &2.16e-04 & 3.31 &0.091 & 0.43 & 0.08\\
  5 &5.64e-04 & 8.63 &0.006 & 0.08 & 0.01\\ \hline    
\end{tabular}
\end{table}

Finally, there is another type of read-only memory that is available for use in CUDA C---Texture memory is cached on chip, so in some situations it will provide higher effective bandwidth by reducing memory requests to off-chip DRAM. More specifically, texture caches are designed for applications where memory access patterns exhibit a great deal of spatial locality. In a computing application, this roughly implies that a thread is likely to read from an address “near” the address that nearby threads read. Now we compare an implementation without using the texture memory on GPU with the implementation using texture memory. From Table~\ref{table_texture}, we can see that, by using texture, one can save about $15\%$ of the computing time if $L$ is large enough. 
\begin{table}[H]
        \centering\caption{Wall times (second) on GPU ($L_{\theta}=L$) with or without using texture memory, where the ``Improvement'' column shows the improvement by using texture-caching.}\label{table_texture}
        \smallskip
        \begin{tabular}{cccc|cccc}\hline
         \multicolumn{4}{c|}{$\mbox{2D}$} & \multicolumn{4}{c} { $\mbox{3D}$}
        \\ \cline{1-8}
      $L$ & No Texture & Texture  & Improvement
 &$L$&  No Texture & Texture  & Improvement
 \\\hline
     9& 1.12e-2 &1.03e-2 &7.62\%  & 5 &7.80e-3 &8.76e-3 &-12.38\% \\
     10& 3.60e-2 &3.10e-2 &13.90\%  & 6 &2.37e-2 &2.36e-2 &0.46\% \\
      11& 1.32e-1 &1.12e-1 &15.10\% & 7 &1.10e-1 &9.92e-2 &9.53\% \\
      12& 5.09e-1 &4.33e-1 &14.98\% & 8 &8.25e-1 &7.26e-1 &12.01\% \\ \hline
       \end{tabular}
\end{table}

\subsection{FMG vs. FFT}\label{subsec:numer_fmg}

Now we compare the FFT method with the geometric multigrid method as fast Poisson solution methods. FFT is a direct solver, and multigrid is an iterative solver. Therefore, making a fair comparison between the two is not an easy task. We set up our comparison in the following way: we tested a sequence of FMG methods, each of which had a different number of pre$-$ and post-relaxation sweeps. Then we compared FFT with the most efficient FMG scheme in order to determine which gives the optimal approximation error. As FFT and FMG require the same amount of data to be transmitted, we only compare the respective kernel times here.

We consider cases with $16$ million unknowns in 2D ($L=12$) and 3D ($L=8$). For the 2D case, from Tables~\ref{table_error_fft_fmg_2d} and \ref{table_kertime_fft_gmg_2D1}, we notice that FMG(1,2) is enough to guarantee the optimal convergence of the approximation error in $L^2(\Omega)$. On the other hand, for the 3D case, we need to use at least FMG(3,3) in order to obtain the optimal convergence rate (see Tables~\ref{table_error_fft_fmg_3d}
and \ref{table_kertime_fft_gmg_3D1}). Moreover, the optimal FMG is $33\%$ and $23\%$ faster than FFT in 2D and 3D, respectively.

\begin{table}[H]
\centering\caption{ Approximation error $\|u-u_h\|$ in 2D}\label{table_error_fft_fmg_2d}
    \begin{tabular}{ccccccc} \hline
    $L$ &FFT      &FMG(1,1) &FMG(1,2)  &FMG(2,2)    &FMG(2,3)&FMG(3,3)  \\ \hline
   9  &1.563e$-$6 &1.001e$-$5  &1.242e$-$6  &1.004e$-$6  &7.028e$-$7  &7.145e$-$7\\ 
   10 &3.914e$-$7 &2.618e$-$6  &3.113e$-$7  &2.518e$-$7  &1.762e$-$7  &1.790e$-$7\\ 
   11 &9.797e$-$8 &6.766e$-$7  &7.791e$-$8  &6.304e$-$8  &4.411e$-$8  &4.479e$-$8\\ 
   12 &2.450e$-$8 &1.735e$-$7  &1.948e$-$8  &1.577e$-$8  &1.103e$-$8  &1.120e$-$8\\ \hline
    \end{tabular}
\end{table}

\begin{table}[H]
\centering\caption{ Kernel time (seconds) of FFT and FMG in 2D}\label{table_kertime_fft_gmg_2D1}
    \begin{tabular}{ccccccc} \hline
    $L$    & FFT       &FMG(1,1)    &FMG(1,2)      &FMG(2,2)    &FMG(2,3)    &FMG(3,3)  \\ \hline
   9          &3.739e$-$3    &3.611e$-$3   &4.260e$-$3   &4.980e$-$3  &5.617e$-$3   &6.348e$-$3\\
   10         &1.102e$-$2    &7.434e$-$3   &8.770e$-$3   &1.008e$-$2  &1.144e$-$2   &1.282e$-$2\\
   11         &4.077e$-$2    &2.203e$-$2   &2.571e$-$2   &2.945e$-$2  &3.317e$-$2   &3.701e$-$2\\
   12         &1.364e$-$1    &7.860e$-$2   &9.167e$-$2   &1.049e$-$1  &1.180e$-$1   &1.310e$-$1\\ \hline
    \end{tabular}
\end{table}

\begin{table}[H]
\centering\caption{ Approximation error $\|u-u_h\|_2$ in 3D}\label{table_error_fft_fmg_3d}
    \begin{tabular}{ccccccc} \hline
    $L$ & FFT    &FMG(1,1) &FMG(1,2)  &FMG(2,2)    &FMG(2,3)    &FMG(3,3)  \\ \hline
    5 & 2.841e$-$4 &6.509e$-$3  &2.733e$-$3  &1.246e$-$3  &7.873e$-$4  &5.296e$-$4\\ 
    6 & 7.100e$-$5  &2.685e$-$3  &9.469e$-$4  &3.930e$-$4  &2.426e$-$4  &1.608e$-$4\\ 
    7 & 1.774e$-$5 &1.032e$-$3  &2.988e$-$4  &1.125e$-$4  &6.751e$-$5  &4.394e$-$5\\ 
    8 & 4.437e$-$6 &3.803e$-$4  &8.880e$-$5  &3.049e$-$5  &1.784e$-$5  &1.145e$-$5\\ \hline
    \end{tabular}
\end{table}

\begin{table}[H]
\centering\caption{ Kernel time (seconds) of FFT and FMG in 3D}\label{table_kertime_fft_gmg_3D1}
    \begin{tabular}{ccccccc} \hline
    $L$    &FFT        &FMG(1,1)     &FMG(1,2)     &FMG(2,2)    &FMG(2,3)    &FMG(3,3)  \\ \hline
   5          &5.102e$-$4   &1.611e$-$3   &1.932e$-$3   &2.382e$-$3  &2.738e$-$3  &3.186e$-$3\\ 
   6          &1.890e$-$3   &3.711e$-$3   &4.474e$-$3   &5.335e$-$3  &6.098e$-$3  &6.986e$-$3\\ 
   7          &5.884e$-$2   &1.342e$-$2   &1.586e$-$2   &1.846e$-$2  &2.094e$-$2  &2.352e$-$2\\ 
   8          &1.893e$-$1   &8.566e$-$2   &1.007e$-$1   &1.155e$-$1  &1.302e$-$1  &1.456e$-$1\\ \hline
    \end{tabular}
\end{table}

%% file: conclusion.tex

\section {Conclusion}\label{sec:conclusion}

In this work, we studied the performance of GMG on CPU--GPU heterogenous computers. Our numerical results suggest that in the best-case scenario the GPU version of GMG can achieve 11 times speed-up in 2D and 10 times speed-up in 3D compared with an efficient OpenMP implementation of multigrid methods on CPUs. When the problem is relatively small we found that the heterogenous algorithm ($0<L_\theta<L$) usually gives the best computational performance. On the other hand, when the problem size is large enough, then it is generally preferable to do the computation on GPUs. We observed lacking of parallelism and frequently visiting coarse levels account for the relatively low floating-point performance of GMG on GPUs. Furthermore, we compared our method with the Fast Fourier Transform in the state-of-the-art cuFFT library. For the test cases with $16$ million unknowns ($L=12$ in 2D and $L=8$ in 3D), we showed that our FMG method is $33\%$ and $23\%$ faster than FFT in 2D and 3D, respectively. Of at least equal importance is that GPU is more cost-effective (in terms of initial cost and daily energy consumption) than modern multicore CPUs for geometric multigrid methods.

%% file: acknowledgements.tex

\section*{Acknowledgements}
The authors would like to thank Dr. Yunrong Zhu from Idaho State University and Dr. Xiaozhe Hu from Penn State University for their helpful comments and suggestions. They are also grateful for the assistance provided by  Mr. Xiaoqiang Yue and  Mr. Zheng Li from Xiangtan University in regard in our  numerical experiments.